\documentclass[12pt,onecolumn ]{IEEEtran}

\usepackage{url}
\usepackage{kbordermatrix}
\usepackage{amsmath,amssymb}
\usepackage{amsthm,thmtools}
\usepackage{amsfonts} 

\usepackage[dvips]{epsfig}      
\usepackage{graphicx}
\usepackage{enumerate}
\usepackage{color}                  

\usepackage{verbatim}
\usepackage{amssymb}
\usepackage{amsbsy}
\usepackage{setspace}
\usepackage{url}

\usepackage{float}

\newcommand{\QED}{$\square$}

\renewcommand{\kbldelim}{(}
\renewcommand{\kbrdelim}{)}

\declaretheorem[name={Example},qed={\lower-0.3ex\hbox{$\square$}} ] {Example}

\declaretheorem[name={Definition}  ] {Definition}
\declaretheorem[name={Theorem}  ] {Theorem}
\declaretheorem[name={Lemma}  ] {Lemma}
\declaretheorem[name={Remark}  ] {Remark}
\declaretheorem[name={Corollary}  ] {Corollary}
\declaretheorem[name={Assumption}  ] {Assumption}

\newcommand {\R}{\mathbb R}

\newcommand {\M}{\mathbb M}

\newcommand {\I}{\mathbb I}

\newcommand{\be}{\begin{equation}}
\newcommand{\ee}{\end{equation}}
\newcommand{\sgn}{\operatorname{{\mathrm sgn}}}
\newcommand{\V}{\mathcal V}

\newcommand{\W}{\mathcal W}

\newcommand*\dif{\mathop{}\!\mathrm{d}}

\begin{document}

\title{ \Large A Generalization of Smillie's Theorem on  Strongly Cooperative Tridiagonal
Systems \thanks{Research  supported in part 
by  research grants from  the Israeli Ministry of Science and Technology,    the Israel Science Foundation, and the 
 US-Israel Binational Science Foundation.  
}} 

\author{Eyal Weiss\thanks{E. Weiss is with the School of Electrical Engineering, Tel Aviv University, Tel-Aviv 69978, Israel.
		E-mail: ‫eyalweiss@mail.tau.ac.il‬}  and	Michael Margaliot\thanks{M. Margaliot (corresponding author) is with the School of Electrical Engineering and the Sagol School of Neuroscience, Tel Aviv University, Tel-Aviv 69978, Israel.
		E-mail: michaelm@post.tau.ac.il }}

\maketitle

\begin{abstract}
Smillie  (1984) proved an interesting  result on the stability of
 nonlinear,  time-invariant, strongly cooperative, and 
tridiagonal dynamical systems. This result has found many applications 
in models from various fields   including biology, ecology, and chemistry.
 Smith (1991) has extended Smillie's  result and proved entrainment 
in the  case where the vector field is time-varying and periodic. 

We use the theory of linear totally nonnegative  differential   systems developed by Schwarz~(1970)
to give a generalization of these two results. This is based on weakening  the
requirement for strong cooperativity  to  cooperativity, and adding an additional observability-type condition. 
\end{abstract}


\section{Introduction}
For two vectors~$a,b\in\R^n$ let~$a\leq b$ denote that~$a_i\leq b_i$ for all~$i$. 
The system~$\dot x=f(x)$ is called 
\emph{cooperative} if  for any two initial conditions~$a\leq b$
the corresponding solutions 
satisfy~$x(t,a)\leq x(t,b)$ for all time~$t\geq 0$. In other words,
the solutions preserve the (partial) ordering~$\leq$ 
 between the initial conditions. The system is
\emph{strongly cooperative} if  for any two  initial conditions~$a\leq b$, with~$a\not =b$,
the corresponding solutions 
satisfy~$x(t,a)< x(t,b)$ for all time~$t> 0$.

Monotone dynamical systems satisfy a similar condition but for a more general ordering~$\leq^K$ that is defined using a suitable cone~$K\subseteq\R^n$.
 Hirsch's quasi-convergence theorem~\cite{hlsmith} shows that monotonicity 
 has far-reaching implications on the asymptotic behavior of the solutions. Roughly speaking, 
it implies that almost all bounded trajectories converge to an equilibrium. 
However, 
monotone  systems can contain complicated dynamics such as chaotic invariant sets (although these
cannot be attractors)~\cite{hirs2}.

Stronger results hold for cooperative systems with an additional structure.
 Let~$\M^+$ denote the set of $n\times n$ matrices that are tridiagonal, and  with positive entries on the super- and sub-diagonals. 
 In an interesting paper, Smillie~\cite{smillie} considered the time-invariant  nonlinear cooperative system:
\be\label{eq:fsys}
					\dot x(t)=f(x(t)),
\ee
where~$x(t) \in \R^n$ and~$f:  \R^n \to \R^n$ is~$n-1$ times differentiable. 
He showed that if the Jacobian~$J(x):=\frac{\partial}{\partial x} f(x) \in \M^+$
for all~$x$,
then \emph{every}   trajectory of~\eqref{eq:fsys} either eventually 
leaves any compact set or
converges to an equilibrium point.\footnote{Note  that Fiedler and  Gedeon~\cite{Fiedler1999} have proved  a similar  result, but using a very different technique. }   
Smillie's result has found many applications (see, e.g.~\cite{RFM_stability,chua_roska_1990,Donnell2009120,HIRSCH1989331}).

Smillie's proof is quite interesting and is based on studying the number of sign variations  
in the vector~$z(t):=\dot x(t)$.
Recall  
	that for a vector~$y\in\R^n$ with no zero entries the number of sign variations
 in~$y$ is
\be\label{eq:sigmdfr}
\sigma(y):=\left|\{i \in \{1,\dots,n-1\} : y_i y_{i+1}<0\} \right | .
\ee
The domain of~$\sigma$    can be extended, via continuity, to the open set
$
\V := \{y\in\R^n:   y_1 \not =0,\; y_n \not=0,\; \text{and if } y_i=0 \text{ for some~$i \in \{2,\dots,n-1\} $ then } y_{i-1}y_{i+1}<0\}.
$
 For example, for~$n=3$ the vector~$y:=\begin{bmatrix}  2 & \varepsilon & -2 \end{bmatrix}' \in \V$ 
and~$\sigma(y)=1$ for \emph{all}~$\varepsilon  \in \R$.

To explain Smillie's  proof,
let~$z:=\dot x$. Then~\eqref{eq:fsys} yields 
\be\label{eq:zdoteqjj}
				\dot z(t)= J(x(t))z(t).
	\ee
This  
 linear time-varying system  
 is sometimes referred to as the variational equation. 

Smillie showed that if~$z(\tau) \not \in \V$ for some time~$\tau$
then
$z(\tau^-),z(\tau^+) \in \V$ and
\[
										\sigma(	z(\tau^+))<\sigma( z(\tau^-)). 
\]
 In other words,~$\sigma(z(t))$ is piecewise constant and whenever it changes its value it can only 
decrease.  
Since~$\sigma$ takes values in~$\{0,\dots,n-1\}$, this implies that~$z(t)  \in \V$
for all~$t\geq 0$ except perhaps for up to~$n-1$
discrete points. By the definition of~$\V$, we conclude that 
there exists a time~$s$ such that~$z_1(t)\not =0 $ (and~$z_n(t)\not = 0 $) for all~$t\geq s$.
Thus, $z_1(t)$ (and~$z_n(t)$) is either eventually  positive or eventually  negative.  
Smillie used this to show that for any~$a$ in the state-space of~\eqref{eq:fsys}
the omega limit set~$\omega(a)$ includes no more than a single point. Hence, every  trajectory either eventfully
 leaves any compact set
or
converges to an equilibrium.

 To explain the basic idea underlying Smillie's analysis of~$\sigma(z(t))$, consider the case~$n=3$. Seeking a contradiction, assume that
	the sign pattern of~$z(t)$ near some time~$t=\tau$ is as follows: 
\begin{center}
\begin{tabular}{c | c c c } 
& $t=\tau^-$ & $t=\tau$ & $t=\tau^+$ \\  \hline 
$z_1(t)$& $+$ & $0$ & $-$ \\  
$z_2(t)$ & $+$ & $+$ & $+$\\    
 $z_3(t)$& $+$ & $+$ & $+$   \\  
\end{tabular} \\
\end{center}
Note that in this case~$\sigma(t):=\sigma(z(t))$ \emph{increases} from~$\sigma(\tau^-)=0$ to~$\sigma(\tau^+)=1$.
However, using~\eqref{eq:zdoteqjj} and  the structure of  the Jacobian yields  
\begin{align*}
							\dot z(\tau)=\begin{bmatrix}  *& + & 0 \\ + & * & +\\ 0& + &* \end{bmatrix}
							\begin{bmatrix} 0\\ + \\ +  \end{bmatrix}
\end{align*}
where~$+ $ [$*$] means a positive [arbitrary]
value, and thus~$\dot z_1(\tau)> 0$, and the 
case described in the table above \emph{cannot} take place. 
	Smillie's  analysis shows rigorously that when~$\sigma(t)$ changes it can only decrease.
	This is  based on direct
	analysis of the~ODEs and is nontrivial 
	  due to the fact that
if an entry~$z_i(t)$ becomes zero at some time~$t=\tau$ (thus perhaps leading to a change in~$\sigma( t  )$ near~$\tau$)
one must consider the possibility that
	higher-order derivatives of~$z_i(t )$ are also zero at~$t=\tau$.
	
	Smith~\cite{periodic_tridi_smith} has extended Smillie's result 
	to time-varying, $T$-periodic cooperative systems under similar assumptions on
		the structure of the Jacobian.
	He showed that every   trajectory 
	eventually leaves any compact set or 
	converges to a periodic trajectory with the  same
	period~$T$.
	This   entrainment property  is important in many natural and artificial systems. 
	For example, biological organisms are often exposed to periodic excitations like the~24h
	solar day, and the periodic cell  division process. Proper functioning often requires entrainment  
	to such excitations~\cite{entrain2011,RFM_entrain}.  
	Epidemics of infectious diseases often correlate with 
	seasonal changes and  the required interventions, such as pulse vaccination, may also need to be periodic~\cite{epidemics_2006}. 
	Traffic flow is often controlled by periodically-varying traffic lights, and in this 
	context entrainment  means that the traffic flow converges to a periodic  pattern with the same period
	as the traffic lights~\cite{2017arXiv171007321M}.

It has been recently shown~\cite{fulltppaper}  that 
  the sign variation diminishing property~(SVDP) of~$z(t)$ underlying these
	results follows as a special case from the seminal (yet largely forgotten)
	work of   Schwarz~\cite{schwarz1970} on  \emph{totally positive differential systems}~(TPDSs).
	These are linear time-varying systems whose transition matrix is totally positive~(TP) for all time. 
	To explain this, recall that a matrix~$A\in\R^{n\times n }$ is called
	totally nonnegative~(TN) if all its minors are nonnegative,
	and totally positive~(TP) if all its minors are positive.
	Such matrices have a rich and beautiful theory~\cite{total_book,pinkus}.
	In particular, they satisfy powerful SVDPs: multiplying a vector by a TN   matrix
	cannot increase
	the number of sign variations in the vector. 
	
	Schwarz~\cite{schwarz1970}  studied 
		the following problem. Fix a time interval~$(a,b)$
	with~$-\infty\leq a<b\leq \infty$.  
	Let~$A(t)$ be a continuous matrix function on~$(a,b)$, and 
	consider the linear matrix differential equation~$\dot \Phi=A \Phi$, 
	with~$\Phi(t_0)=I$. When will~$\Phi(t)$ be TN [TP] for every pair~$(t_0,t)$ with~$a<t_0<  t<b$? 
	A system that satisfies this property is called a
	totally nonnegative [positive] differential system. 
	Schwarz also described the implications of TNDS [TPDS] on the sign variations of the vector
	solution~$z(t)$
	of~$\dot z=A z$.
	
	As shown in~\cite{fulltppaper}, these results are closely related 
	to the theorems  of Smillie, Smith, and others. 
	Indeed, the assumptions of Smillie on the Jacobian   imply that~\eqref{eq:zdoteqjj} is~TPDS.
	However, the work of 
	Schwarz has been largely forgotten and its implications to the analysis of \emph{nonlinear} 
	dynamical systems have    been overlooked.

Here we   generalize one of the results of Schwarz on TNDSs, and then use this to generalize
the results of Smillie and Smith. Roughly speaking, this generalization
is based on weakening  the requirement for a  triadogonal 
Jacobian with positive entries on the super- and sub-diagonal 
 to requiring a  tridiagonal Jacobian with \emph{nonnegative} entries on the super- and sub-diagonals,
  and adding a suitable  observability-type condition.

The next section briefly 
reviews totally positive and totally nonnegative linear differential systems and their
properties. Section~\ref{sec:main} describes our main results, and the final section concludes. 

\section{Preliminaries}
		
		We  begin with reviewing    definitions and results from  the 
		theory of totally nonnegative and totally positive 
		matrices that will be used later on. We consider only square and real matrices, 
		as this is the case that is relevant for our applications. For more information and proofs
		see the     monographs~\cite{total_book,pinkus} and the survey paper~\cite{ando_survey}.
		Unfortunately, this field suffers from nonuniform terminology. We follow the more modern terminology used in~\cite{total_book}.

	 	\begin{Definition}\label{def:tntp}
	 		A matrix~$A\in\R^{n\times n}$
	 		is called totally nonnegative [totally positive]
	 		if the determinant of \emph{every} square submatrix is nonnegative [positive]. 
	 	\end{Definition}
		
		In particular, if~$A$  is~TN [TP] then every
		entry of~$A$   is nonnegative [positive].

		Some matrices with a special structure are known to be~TN.	
		We  review two such examples. The first is important 
				in  proving the SVDP of~TN matrices. 
				The second example  is closely related to Smillie's results.

		\begin{Example}
		Let~$E_{i,j}\in\R^{n\times n}$ denote the  matrix with all entries zero, except for
		entry~$(i,j)$  that is one.
		For~$p \in \R$  and~$i\in\{2,\dots,n\}$, let
							\begin{align}\label{eq:eb}
							L_i(p)&:=I+p E_{i,i-1},\\
							U_i(p)&:=I+p E_{i-1,i}.\nonumber
							\end{align}
							Matrices in this form are called 
		\emph{elementary bidiagonal}~(EB) matrices. If the identity matrix~$I$
		in~\eqref{eq:eb} is replaced by a diagonal matrix~$D$
		 then the matrices are called
		\emph{generalized elementary bidiagonal}~(GEB).
    It is straightforward to see  that~EB matrices 
		are TN when~$p\geq 0$, and that~GEB  matrices  are TN when~$p\geq 0$ and the diagonal matrix~$D$ is componentwise nonnegative.
					\end{Example}

\begin{Example}\label{exa:trid}
		Consider the    tridiagonal    matrix
		\be\label{eq:trida}
							A=\begin{bmatrix} 
																	a_1 & b_1 & 0 & \dots & 0											 \\
																	c_1 & a_2 & \ddots & \dots & \vdots											 \\
																	0 & \ddots & \ddots & \dots & \vdots											 \\
																	\vdots & \ddots & \ddots & \dots & b_{n-1}											 \\
																	0 & \dots & \dots & c_{n-1} & a_{n}											 
									\end{bmatrix} 
		\ee
		where~$ b_i,c_i\geq0$  for all~$i$.
		In this case, the dominance condition
		\be\label{eq:domi}
		a_i \geq b_i+c_{i-1} \quad \text{for all } i\in\{1,\dots,n\},
		\ee
		with~$c_0:=0$ and~$b_n:=0$, guarantees that~$A$ is TN~\cite[Ch.~0]{total_book}.
		\end{Example}

		An important subclass of TN matrices that are ``close'' to TP matrices  
		are the \emph{oscillatory matrices}  studied in the pioneering work
		of Gantmacher and Krein~\cite{gk_book}. 
		A matrix~$A\in\R^{n\times n}$ is called \emph{oscillatory} if   $A$ is~TN and there exists an integer~$k>0$ such that~$A^k$ is~TP. It is well-known that a~TN matrix~$A$ is oscillatory if and only if it 
		is non-singular and irreducible~\cite[Ch.~2]{total_book}, and that in this case~$A^{n-1}$ is~TP. 
		\begin{Example}
		Consider the  matrix
		$
		A=\begin{bmatrix}   1&\varepsilon &0 \\ \varepsilon &1&\varepsilon\\0 &\varepsilon &1   \end{bmatrix},
		$
		with~$\varepsilon \in (0,1/2)$. This matrix is   non-singular  (as $\det(A)=1-2 \varepsilon^2$),
		TN (by the result in Example~\ref{exa:trid}),   and irreducible,
		so it is an oscillatory matrix. Here
	$
		A^{n-1}=A^2=\begin{bmatrix}   1+\varepsilon^2 &2 \varepsilon &\varepsilon^2 \\ 
		2 \varepsilon& 1+2\varepsilon^2 &2 \varepsilon \\\varepsilon^2 &2\varepsilon & 1+\varepsilon^2   \end{bmatrix},
		$
		and it is straightforward to verify that this matrix is indeed~TP. 
		\end{Example}

 More generally,   the matrix~$A$ in~\eqref{eq:trida} with~$b_i,c_i>0$ and the dominance condition~\eqref{eq:domi}
	is~TN and irreducible. If it is also non-singular then it  is oscillatory.	

		An important property, that will be used throughout, is that
			the  product of two~TN [TP] matrices is a~TN [TP]  matrix. 
			This follows immediately from Definition~\ref{def:tntp}
			and the
		Cauchy-Binet formula for the minors of the product of two matrices~\cite[Ch.~0]{matrx_ana}. 
	 	
		When using TN matrices to study dynamical systems, it is important to bear in mind that  in general 
		coordinate transformations do not preserve~TN. 
		An important exception, however, is positive diagonal scaling:
		if~$D $ is a diagonal matrix with positive
		entries on the diagonal then multiplying a matrix~$A$ by~$D$ either on the left or
right   changes the sign of no minor, so in particular~$DAD^{-1}$ is TN [TP] if and only if~$A$ is~TN [TP].
		  
	\subsection{Sign variation diminishing property}
		As noted above~$\sigma(y)$ is not well-defined for all~$y$. We recall 
		two definitions of the number of sign variations that are well-defined for all~$y\in\R^n$.
 	Let~$s^-(y)$ denote the number of sign variations
	in the vector~$y$ after deleting all its zero entries,
	and let~$s^+(y)$ denote the \emph{maximal} possible number of sign variations
	in~$y$ after each zero entry is replaced by either~$+1$ or~$-1$. 
	Note that~$s^-(y) \leq s^+(y)$ for all~$y\in\R^n$. 
	For example, for~$y=\begin{bmatrix} 2& 0 &1 &-2 &0 &2.3 \end{bmatrix}'$, $s^-(y)=2$ and~$s^+(y)=4$. 
	 Let~$
			\W:=\{y\in\R^n:s^-(y)=s^+(y)\}.
$
Note   that  if~$y\in \W$ then~$y$ cannot have two adjacent zero coordinates.
An immediate yet important observation is that~$\W=\V$.

		Let~$A$ be a TN  EB matrix. 
		Pick~$x\in \R^n$, and let~$y:=Ax$. Then there exists at most one index~$i$ such  
		  that~$\sgn(y_i) \not = \sgn(x_i)$, and
			either~$y_i=x_i+p x_{i-1}$ or~$y_i=x_i+p x_{i+1}$, and since~$p\geq 0$
		the sign can change only in the ``direction'' of~$x_{i-1}$ or~$x_{i+1}$. 
		In either case, neither~$s^-$ or~$s^+$ may increase. 
		We  conclude that if~$A$ is TN  EB  then 
	$		s^-(Ax) \leq s^-(x),
		$
		and
		\be\label{eq:splushcnage}
		s^+(Ax) \leq s^+(x) .
		\ee
		A similar argument shows that if~$A$ is TN GEB then~$		s^-(Ax) \leq s^-(x)$.
		However,~\eqref{eq:splushcnage} does not hold in general for a TN GEB matrix~$A$. 
		For example,~$A=0$
		 is   TN GEB and clearly~$s^{+}(Ax)=s^{+}(0)=n-1$ may be larger than~$s^{+}(x)$.
 		
			This  SVDP can be extended to all~TN matrices using   
			a fundamental   decomposition result  that states that 
		any TN matrix  can be expressed as a product of TN GEB matrices~\cite[Ch.~2]{total_book}.
	\begin{Theorem}\cite[Ch.~4]{total_book}  \label{Pthm:svdpo}  
If~$A\in\R^{n\times n}$ is TN then 
				\be\label{eq:smtn}
				s^{-}(Ax)\leq s^{-}(x),\text{ for all } x\in\R^n.
				\ee
				If $A$ is TN and  nonsingular then 
					\be\label{eq:spod} 
						s^{+}(Ax)\leq s^{+}(x), \text{ for all } x\in\R^n.
						\ee
							If~$A $ is TP then 
					\be\label{eq:strong}
									s^+(Ax)\leq s^-(x),\text{ for  all }  x\in \R^n \setminus\{0\}.
					\ee
					If~$A$ is TN and nonsingular then~$s^+(Ax)\leq s^-(x)$  holds  
					if either~$x$ has no zero entries or~$Ax$ has no zero entries. 

				\end{Theorem}

					At this point we can already reinterpret Simillie's result  using the SVDP of TP matrices. 
					To explain this, consider for simplicity the system~$\dot z= J z$, with~$J$ a  \emph{constant}
					 matrix and~$J\in\M^+$. 
					Then clearly there exists~$\bar \varepsilon>0$ sufficiently small such that 
					the matrix~$\exp(J\varepsilon)=I+\varepsilon J+o(\varepsilon)$
					  is nonsingular, TN (by the result in Example~\ref{exa:trid}),
						and irreducible for all~$\varepsilon  \in (0,\bar \varepsilon)$. Thus,~$(\exp(J\varepsilon))^{n-1}$
					is TP, implying that~$ \exp(J\varepsilon)   $ is in fact TP for all~$\varepsilon>0$ sufficiently small. 
					Since~$z(\varepsilon) =\exp(J\varepsilon) z(0)$, we   conclude from Theorem~\ref{Pthm:svdpo} that~$s^+(z(\varepsilon))\leq s^-(z(0))$ for all~$z(0)\not =0$.

					This suggests the following  question:
					 when is the transition matrix associated with~$\dot x=A x$ TN [or TP]?
					This is exactly the question   addressed by Schwarz in~\cite{schwarz1970}. He introduced the following definitions.\footnote{We slightly modify the original definitions in~\cite{schwarz1970} to make them compatible with   more
					modern terminology.}
					
\begin{Definition}
				Consider the matrix differential  equation~$\dot \Phi=A \Phi$, $\Phi(t_0)=I$, where~$A(t)$ is a continuous 
				matrix on the time interval~$t \in (a,b)$. 
					The system is called a totally nonnegative differential system~(TNDS) if~$\Phi(t)$ is~TN for any pair~$(t_0,t)$ with~$a<t_0\leq t   <b$. It is called 
						  a totally positive differential system~(TPDS) if~$\Phi(t)$ is~TP for any pair~$(t_0,t)$ with~$a<t_0< t <  b$.
\end{Definition}

For the case where~$A(t)$ is continuous in~$t$, Schwarz derived a necessary and sufficient condition for a
 system to be~TNDS or~TPDS. 
						 Let~$\M \subset \R^{n\times n}$ denote the set of tridiagonal matrices
with nonnegative entries on the sub- and super-diagonals.

					\begin{Theorem}\cite{schwarz1970}
					Consider the matrix differential  system~$\dot \Phi=A \Phi$, $\Phi(t_0)=I$, where~$A(t)$ is a 
					continuous matrix for~$t\in(a,b)$. 
					The system is TNDS iff~$A(t) \in \M$ for all~$t\in(a,b)$. It is 
					  TPDS iff~$A(t) \in \M$ for all~$t\in(a,b)$,  and every entry on the super- and sub-diagonals of~$A(t)$ is not zero on a time interval.  
					\end{Theorem}
					
					Schwarz also analyzed the implications of TNDS or TPDS to the sign variations in the vector solution of~$\dot z(t)=A(t) z(t)$. 
					
					The next section describes our main results. These are based on weakening 
					the requirement~$J(t) \in \M^+$ for all~$t$ to~$J(t) \in \M$ for all~$t$ and adding an observability-type condition.
					
\section{Main results}\label{sec:main}
Our first result provides a  bound on the number of isolated zeros of~$z_1(t)$ and~$z_n(t)$ 
 in the system~$\dot z = Az$ that  is~TNDS (but not necessarily~TPDS).  
From here on we assume  a more genral case than in Schwarz~\cite{schwarz1970}, namely, that
\be\label{eq:sas}
\text{ }  A:(a,b) \to \R^{n\times n} \text{ } 
\parbox{20em}{ is a  matrix of locally (essentially) \\bounded  measurable 
	functions.}
\ee
 Recall that~\eqref{eq:sas} implies that 
\be\label{eq:phidd}
\dot\Phi=A\Phi, \quad \Phi(t_0)=I,
\ee
 admits a unique, locally absolutely continuous,    invertible
solution for   all~$t \in(a,b)$ 
(see, e.g.,~\cite[Appendix~C]{sontag_book}).




	\begin{Theorem}\label{thm:exp_is_tnn}
Consider the time-varying linear system:
\be\label{eq:linti}
\dot z(t)=A(t)z(t), 
\ee
with~$A(t)$ satisfying~\eqref{eq:sas} 
and suppose that
\be \label{eq:phist}
\Phi(t,t_0) \text{ is TN for all } a<  t_0 \leq t<b   . 
\ee 
Let~$p_i$ denote the number of \emph{isolated} zeros of~$z_i(t)$ on~$(a,b)$.
Then~$\max\{p_1,p_n\}\leq n-1$. 
\end{Theorem}

{\sl Proof of Thm.~\ref{thm:exp_is_tnn}.} 
We prove the assertion for~$z_1(t)$ (the proof for~$z_n(t)$ is very similar).
Since~$\Phi(t)$ is~TN and invertible for all~$t\geq t_0$,
Thm.~\ref{Pthm:svdpo}  implies that~$s^+(z(t))$ is non-increasing with~$t$. 	
Schwarz~\cite{schwarz1970} has shown that if~$\dot z=Az$ is a TNDS on~$(a,b)$ and  there
exist
  times~$r,q$ with~$a<r<q<b$ such that
$
z_1(r)=0$ and $ z_1(q)\not =0$
then
		$	s^+(z(q)) < s^+(z(r)) $. 
Since~$s^+$ takes values in~$\{0,1,\dots,n-1\}$, we conclude that~$z_1(t)$ cannot have more than~$n-1$
isolated zeros on~$(a,b)$.~$\square$


Note that Thm.~\ref{thm:exp_is_tnn} only bounds the number of  \emph{isolated} zeros of~$z_1(t)$ and~$z_n(t)$. In particular, it does not rule out the possibility that~$z_1(t) $ or~$z_n(t)$
are zero on a time interval. 
\begin{Example}
Consider the case~$n=2$ and the constant  matrix
$
			A = \begin{bmatrix} a_{11} & 0 \\ 1   &a_{22} \end{bmatrix}.
$
  Then~$A\in\M$, so~\eqref{eq:linti} is~TNDS.
The solution of~$\dot \Phi =A \Phi$, $\Phi(0)=I$, is
\[
			   \exp(At)=\begin{bmatrix}\exp( a_{11}t)   & 0 \\
			p(t)  &\exp(a_{22}t) \end{bmatrix} ,
\]
where~$p(t):=\frac{\exp(a_{11}t)-\exp(a_{22}t) }{a_{11}-a_{22}} $ if~$a_{11}\not= a_{22}$,
and~$p(t):=t \exp(a_{11}t) $, otherwise. Note that in both cases~$p(t)>0$ for all~$t>0$, so~$\exp(At)$
is~TN for all~$t\geq 0$.
From this we conclude that the solution of~$\dot z=Az$ 
is~$z(t)= \begin{bmatrix} 0 &  \exp(a_{22}t)z_2(0) \end{bmatrix}'$
if~$z_1(0)=0$, 
and~$z(t)= \begin{bmatrix} \exp(a_{11}t)z_1(0)& p(t)z_1(0)+\exp(a_{22}t)z_2(0)\end{bmatrix}'$
if~$z_1(0)\not =0$.


In the first case,~$z_1(t)$ is zero for all~$t\geq 0$ and~$z_2(t)$ is either zero for all~$t$
or has no zeros. In the second case, $z_1(t)$ has no zeros, and~$z_2(t)$ has no more than a single isolated zero. 
\end{Example}

\subsection{Applications to stability analysis}\label{sec:app}
We are now ready to give a generalization of Smillie's Theorem. 
In fact, we provide a generalization to a result of Smith on the stability of tridiagonal
cooperative systems 
with a time-varying periodic vector field, and then specialize to the time-invariant case.

Consider  the  nonlinear  time-varying dynamical system 
\be\label{eq:gfd}
\dot  x(t)=f(t,  x(t)),
\ee
whose  trajectories   evolve  on an invariant
 set~$\Omega \subset \R^n$, that is, for any~$x_0\in\Omega$ and any~$t_0\geq 0$
a unique solution~$x(t,t_0,x_0)$ exists and satisfies~$x(t,t_0,x_0)\in\Omega$ for all~$t\geq t_0$. From here on we take~$t_0=0$.  We assume that the state-space~$\Omega$ is convex and compact.   
We also assume  that the 
Jacobian~$J(t,x)=\frac{\partial}{\partial  x} f(t,x)$
exists for all~$t\geq 0$ and~$x\in \Omega$.

We consider the case where~$f$ is~$T$-periodic, that is, 
\[
f(t,z)=f(t+T,z), \text { for all } t\geq 0,\; z\in \Omega. 
\]
Note that in the particular case where~$f$ is time-invariant this property  holds for all~$T$.
We make two assumptions.

					\begin{Assumption}\label{assump:atmat}
For any~$t\geq 0 $ and any line~$\gamma:[0,1]\to \Omega$   
  the matrix
\be\label{eq:atpo}
A(t):=\int_0^1  J(t,\gamma(r))     \dif r
\ee
 is well-defined,  locally (essentially) bounded,  measurable, and satisfies the conditions for~TNDS, i.e.~$A(t) \in \M$ for almost all~$t \in(a,b)$. 
\end{Assumption}
	Note that since~$f$ is~$T$-periodic so is~$A(t)$.

					\begin{Assumption}\label{assume:sec}
For any solution of~$\dot z(t) =A (t) z(t)$ that is not the trivial
					solution either~$z_1(t)$ or~$z_n(t)$ has only isolated zeros.
\end{Assumption}

\begin{Theorem}\label{thm:application}
If Assumptions~\ref{assump:atmat} and~\ref{assume:sec}  hold  then 
every solution of~\eqref{eq:gfd} converges to a periodic trajectory with period~$T$. 
\end{Theorem}

\begin{Remark}
Thm.~\ref{thm:application} is a generalization of a result of Smith~\cite{periodic_tridi_smith} who derived the same result
under stronger   conditions, namely, that~$A(t)$ satisfies the necessary conditions for~\emph{TPDS} (so in particular Assumption~\ref{assump:atmat} holds).
TPDS also means  that the solution of~$\dot z=Az$ satisfies~$z(t) \in \V$ 
for all~$t\geq 0$ except perhaps for up to~$n-1$ discrete time points~\cite{schwarz1970}.
  By the definition of~$\V$, 
this means that Assumption~\ref{assume:sec} also 
holds. Thus, the result of Smith is a special case of
Thm.~\ref{thm:application}.
\end{Remark}

If the vector field has the form~$f(t,x)=f(x,u)$, with~$u(t)$  $T$-periodic, one may view~$u$ as a periodic input
that excites the system.
Thm.~\ref{thm:application} then implies that the system \emph{entrains} to the excitation, as
every solution also converges to a periodic solution with the same period as the excitation. 

We will prove Thm.~\ref{thm:application} in the case where all the zeros of~$z_1(t)$ are isolated (the proof in the case where the zeros of~$z_n(t)$ are isolated is very similar).
The proof  is similar to the proof in Smith~\cite{periodic_tridi_smith} (albeit under our weaker assumptions), but we include it for the sake of completeness. 
We require the following result.

\begin{Lemma}\label{lem:ztr}
Pick~$a,b \in \Omega$, with~$a\not = b$,  and consider the solutions~$x(t,a)$, $x(t,b)$  
  of~\eqref{eq:gfd}. Then there exists a time~$s\geq 0$ such that for all~$t\geq s$
	either~$x_1(t,a)>x_1(t,b)$ or~$x_1(t,a)<x_1(t,b)$. 
	\end{Lemma}
	 
	{\sl Proof of Lemma~\ref{lem:ztr}.}
Let
	\[
	\gamma(t,r):=rx(t,a)+(1-r)x(t,b),\quad r\in[0,1],
	\]
	denote the line between the two solutions at time~$t$.
	Since~$\Omega$ is convex,~$\gamma (t,r) \in \Omega$ for all~$t\geq 0$ and all~$r\in[0,1]$. 
Let~$z(t):=x(t,a)-x(t,b)$. Then
\begin{align*}
   \dot z(t)&=f(t,x(t,a))-f(t,x(t,b))\\
	       &=\int_0^1  \frac{d }{d r}   f(t,\gamma(r))     \dif r\\
					&=A(t)z(t),
\end{align*}
with~$A(t)$ defined  in~\eqref{eq:atpo}.
 By Assumption~\ref{assump:atmat}, this is a TNDS. Hence, according to  
 Thm.~\ref{thm:exp_is_tnn}, $z_1(t)$ has no more than~$n-1$ isolated zeros. 
 Combining this with Assumption~\ref{assume:sec}, which states that 
 $z_1(t)$ has only isolated zeros, we conclude  that
 there exists~$s\geq 0$ such that~$z_1(t)\not = 0  $ for all~$t\geq s$ and this completes the proof.~\hfill{\QED}
 
We can now prove  Thm.~\ref{thm:application}.   
If for some~$a\in \Omega$ the solution~$x(t,a)$ is $T$-periodic then there is nothing to prove. Thus,
suppose that for some~$a\in \Omega$ the solution~$x(t,a)$ of~\eqref{eq:gfd}  is  \emph{not}~$T$-periodic. 
  Then the~$T$-periodicity of the vector field implies that~$ x(t+T,a)$ is another solution of~\eqref{eq:gfd} that is different from~$x(t,a)$.
	 Lemma~\ref{lem:ztr} implies that
	there exists 
  an integer~$m\geq 0$ such that~$x_1(k T ,a)-x_1((k+1)T,a) \not =0$ for all~$k\geq m$. Without loss of generality,
	  assume that 
\be\label{eq:xms}
					x_1(k T,a )-x_1((k+1)T,a) >0 \text{ for all } k\geq m. 
\ee
 
   Define the Poincar\'e
 map~$P_T:\Omega\to \Omega$ by~$
				P_T(y):=x(T, y).
$
Then~$P_T$ is continuous, and for any integer~$k\geq 1$
 the $k$-times composition of~$P_T$ satisfies~$P_T^k(y)=x(kT,y)$. 
 The  omega limit set $\omega_T:\Omega\to \Omega$ is defined by
 $
 \omega_T(y): = \{ z\in \Omega:  \text{ there exists a sequence } n_1,n_2, \dots  
 \text{ with } n_k \to \infty 
 \text{ and }\lim_{k\to \infty} P_T^{n_k}(y)= z \} . 
 $  
It is well-known that~$\omega_T(y)\not =\emptyset$, $  x(kT,y)  \to  \omega_T(y)$, 
 and $\omega_T(y)$ is  invariant under~$P_T$, that is,~$P_T( \omega_T(y))  =\omega_T(y)$. 
In particular, if~$\omega_T(y)=\{q\}$ then~$P_T(q)=q$, that is, 
 the solution emanating from~$q$ is $T$-periodic. 
To prove the theorem we need to show that~$\omega_T(a)$ is a singleton. 
Assume that this is not the case. Then there exist~$p, q 
  \in \omega_T(a)$ with~$p\not = q$. 
By the definition of~$\omega_T(a)$, there exist  integer sequences~$n_k\to \infty$ and~$m_k \to \infty$
such that
$
\lim_{k\to \infty}x(n_k T ,a)=p$ and~$\lim_{ k \to \infty}x(m_k T,a)=q. $
Combining this with the monotonicity condition~\eqref{eq:xms}  yields~$p_1=q_1$. 
We conclude that all  points in~$ \omega_T(a)$ have the same first coordinate. 
Now consider the solutions emanating from~$p$ and from~$q$ at time zero, that is,~$x(t,p)$ and~$x(t,q)$.  
We know that  there exists an integer~$m\geq 0$ such that, say, 
\[
						x_1(kT,p)-x_1(kT,q)>0 \text{ for all } k\geq m.
\]
But since~$p,q\in \omega_T(a)$, $x(kT,p),x(kT,q)\in \omega_T(a)$ for all~$k$, and thus we already know that they have the same first coordinate, that is,  $x_1(kT,p)=x_1(kT,q)$. This contradiction 
completes the proof  of Thm.~\ref{thm:application}.~\hfill{\QED}


The time-invariant nonlinear   dynamical system:
\be\label{eq:gfdti}
\dot  x(t)=f( x(t))
\ee
is~$T$ periodic for all~$T>0$, so Thm.~\ref{thm:application} yields the following result.

\begin{Corollary}
Suppose that the solutions of~\eqref{eq:gfdti}  evolve  on an invariant  compact and convex set~$\Omega \subset \R^n$,
that  
  the matrix~$J(x):=\frac{\partial}{\partial x}f(x)\in\M $   
   for all~$x\in \Omega$, and that
	Assumption~\ref{assume:sec} holds. 
	Then for every~$x_0\in\Omega$ the solution~$x(t,x_0)$
  converges to an equilibrium point. 
\end{Corollary}

This is  a generalization of Smillie's theorem~\cite{smillie}.
Indeed, Smillie assumed that~$J(x)\in \M^+$ for all~$x\in\Omega$.
Note that~$J(x)\in \M$   means   that~$J(x) $ 
is tridiagonal and  Metzler, so~\eqref{eq:gfd}
  is  a tridigonal   cooperative  
system in the sense of~Hirsch (see~\cite{hlsmith}).

 To establish that  Assumption~\ref{assume:sec}  indeed holds 
one can use an observability-type test. 
Indeed, consider the  general time-varying, nonlinear, cooperative, tridiagonal  dynamical system:
\begin{align*}
					\dot x_1 &= f_1(t,x_1,x_2),\\
					\dot x_2 &= f_2(t,x_1,x_2,x_3),\\
					\vdots\\
					\dot x_{n-1} &= f_{n-1} (t, x_{n-2},x_{n-1},x_n    ), \\
					\dot x_{n} &= f_{n} ( t, x_{n-1},x_n    ).
\end{align*}
Then~$z:=\dot x$ satisfies
\be\label{eq:zjj}
\dot z=J(t, x(t))z,
\ee
 with
\[
				J=\begin{bmatrix} 
				\frac{\partial f_1}{\partial x_1}   &  \frac{\partial f_1}{\partial x_2}  &0&0&\dots&0\\ 
				\frac{\partial f_2}{\partial x_1}    &  \frac{\partial f_2}{\partial x_2}   &\frac{\partial f_2}{\partial x_3}  &0&\dots&0\\ &&\vdots\\
				0&0&\dots&\frac{\partial f_{n-1}}{\partial x_{n-2}}   &  \frac{\partial f_{n-1}}{\partial x_{n-1}}   &\frac{\partial f_{n-1}}{\partial x_n} \\ 
							0&0&\dots&0&	\frac{\partial f_n}{\partial x_{n-1}}   &  \frac{\partial f_n}{\partial x_n}  \\ 
				\end{bmatrix}.
\]

Suppose that
\be\label{eq:fifi}
\frac{\partial f_i}{\partial x_{i+1}}(t,z) >0 \text{ for all } t\geq 0, z\in \Omega,   i=1,\dots,n-1.
\ee
In other words, the entries on the super-diagonal of~$J$ are positive.

Suppose that~$z_1(t) =  0$ on a time interval~$\I\in(a,b)$. 
Then
\begin{align*}
0 &= \dot z_1  =
   \frac{\partial f_1}{\partial x_2}  z_2,
\end{align*}
and using~\eqref{eq:fifi} implies that~$z_2(t) = 0$ on~$\I$. Thus,
\begin{align*}
0 &= \dot z_2 = 
    \frac{\partial f_2}{\partial x_3}  z_3.
\end{align*}
Proceeding in this manner, we conclude that~$z_1(t)$ is zero on~$\I$ only if~$z(t)$ is the trivial solution, i.e.~\eqref{eq:fifi} implies that
Assumption~\ref{assume:sec}  indeed holds.

Note that unlike in the results by Smillie and Smith we do not require   any entry on the
sub-diagonal of~$J$ to be positive. 
 This demonstrates the fact that our results are more general, as 
if the entries on the sub-diagonal are only nonnegative
 then~\eqref{eq:zjj}
 is~TNDS but not necessarily~TPDS.

More generally,
to check whether Assumption~\ref{assume:sec} holds,  
we can associate with~\eqref{eq:zjj} an output~$y:=z_1$,  
i.e. $y=c'z$, with~$c:=\begin{bmatrix} 1&0&0 &\dots & 0 \end{bmatrix}'$.
Clearly, if~$z_1(t)$ is zero on a time interval~$\I$ and~$z(t)$ is not the trivial solution 
then~\eqref{eq:zjj} is not observable.
Thus, observability of~\eqref{eq:zjj}   on any time interval implies that 
Assumption~\ref{assume:sec} holds. 
One can then apply well-known results guaranteeing
the observability of time-varying linear systems (see,
 e.g.,  \cite{observ_time_varying},\cite[Chapter~6]{sontag_book})
to establish that 
Assumption~\ref{assume:sec} holds. 
Of course, a  similar approach can be used to establish that~$z_n(t)$ is not zero on a time interval.


\section{Conclusions}
Entrainment is an important asymptotic property of dynamical systems and has many applications.
In this paper we considered two interesting results derived
 by Smillie~\cite{smillie} and Smith~\cite{periodic_tridi_smith}, 
that guarantee  stability and entrainment for nonlinear cooperative tridiagonal dynamical systems.
Building upon the theory of  TNDSs developed  by Schwarz~\cite{schwarz1970},
we were able to generalize these stability and entrainment results under 
   a weaker condition, namely that  the Jacobian is tridiagonal but 
may have  nonnegative  (rather than positive) entries on its 
super- and sub-diagonals, along with a suitable observability-type condition.

As a topic for further research, we believe that combining the~TNDS
framework with an additional observability-type condition may be used to generalize other results
that are based on using 
	the number of sign variations in the vector of derivatives
	as a discrete-valued Lyapunov function  (see, e.g.,~\cite{Fusco1990,smith_sign_changes}).

\bibliographystyle{IEEEtranS}
\bibliography{tpds_bib}

\end{document}